# A CENTRAL LIMIT THEOREM VIA DIFFERENTIAL EQUATIONS

By Taral Guldahl Seierstad[1]

*University of Oslo*

In a paper from 1995, Wormald gave general criteria for certain parameters in a family of discrete random processes to converge to the solution of a system of differential equations. Based on this method, we show that if some further conditions are satisfied, the parameters converge to a multivariate normal distribution.

**1. Main theorem.** In this paper, we consider parameters defined on random discrete processes. When the parameters change by only a small amount from one state in the process to the next, one often finds that the parameters satisfy a law of large numbers, that is, the parameters are sharply concentrated around certain values. Wormald [6] gives some general criteria which ensure that given parameters converge in probability to the solution of a system of differential equations.

In fact, such parameters often satisfy not only a law of large numbers, but also a central limit theorem. Based on the differential equation method described in [6] and a martingale central limit theorem due to McLeish [3], we show that when certain general criteria are satisfied, a set of parameters defined on a family of discrete random processes converges to a multivariate normal distribution.

As examples of processes to which this method can be applied, we consider in Sections 4 and 5 two random graph processes. In both processes, the initial state is an empty graph on $n$ vertices, and edges are added one by one according to a random procedure.

Consider a sequence $(\Omega_n, \mathcal{F}_n, P_n)$ of probability spaces. Let $m_n$ be a sequence of numbers such that $m_n = O(n)$, and suppose that for each $n$ a filtration $\mathcal{F}_{n,0} \subseteq \mathcal{F}_{n,1} \subseteq \cdots \subseteq \mathcal{F}_{n,m_n} \subseteq \mathcal{F}_n$ is given. Let $\{\mathbf{X}_{n,m}; m = 0, 1, \ldots, m_n\}$ be a sequence of random vectors in $\mathbb{R}^q$, for some $q \geq 1$, such that $\mathbf{X}_{n,m}$ is

Received October 2007; revised May 2008.

[1]Supported by the Research Council of Norway Grant 170620/V30.

measurable with respect to $\mathcal{F}_{n,m}$ for $0 \le m \le m_n$. The $k$th entry in $\mathbf{X}_{n,m}$ is denoted by $X_{n,m,k}$. For $1 \le m \le m_n$ and $1 \le k \le q$, we define $\Delta\mathbf{X}_{n,m} = \mathbf{X}_{n,m} - \mathbf{X}_{n,m-1}$ and $\Delta X_{n,m,k} = X_{n,m,k} - X_{n,m-1,k}$. If $\mathbf{v}$ is a vector, we let $\mathbf{v}'$ be the transpose of $\mathbf{v}$. If $\mathbf{v}$ is a column vector, we use the notation $\mathbf{v}^2$ to mean $\mathbf{v}\mathbf{v}'$. Thus, if $\mathbf{v}$ is a $q$-dimensional vector, $\mathbf{v}^2$ is a $q \times q$-matrix. We use the norm $\|\mathbf{v}\| = \|\mathbf{v}\|_\infty$. When we use the notation $O(\cdot)$ and $o(\cdot)$, we mean that the bounds hold as $n \to \infty$, unless stated otherwise; if the notation is used to refer to matrices or vectors, the bounds are meant to apply to every entry in the matrix or vector.

If $D \subset \mathbb{R}^q$, we define the stopping time $H_D = H_D(\mathbf{X}_{n,m})$ to be the minimum $m$ such that $n^{-1}\mathbf{X}_{n,m} \notin D$.

The object of this paper is to find criteria which ensure that $\mathbf{X}_{n,\lfloor tn \rfloor}$ converges to a multivariate normal distribution, whose mean and covariance matrix are continuous functions of $t$ and can be obtained by solving certain differential equations. The mean is obtained by applying the following theorem due to Wormald.

Theorem 1 (Theorem 5.1 in [7]). *Assume that there is a constant $C_0$ such that $X_{n,m,k} \le C_0 n$ a.s. for all $n$, $0 \le m \le m_n$ and $1 \le k \le q$. Let $f_k : \mathbb{R}^q \to \mathbb{R}$, $1 \le k \le q$, be functions and assume that the following three conditions hold, where $D$ is some bounded connected open set containing the closure of*

$$\{(z_1, \dots, z_q) : \mathbb{P}[X_{n,0,k} = z_k n, 1 \le k \le q] \ne 0 \text{ for some } n\}.$$

(i) *For some function $\beta = \beta(n) \ge 1$, $\|\Delta\mathbf{X}_{n,m}\| \le \beta$, a.s. for $1 \le m < H_D$.*

(ii) *For some function $\lambda_1 = \lambda_1(n) = o(1)$ and all $k$ with $1 \le k \le q$,*

$$|\mathbb{E}[\Delta X_{n,m,k} \mid \mathcal{F}_{m-1}] - f_k(n^{-1}X_{n,m-1,1}, \dots, n^{-1}X_{n,m-1,q})| \le \lambda_1$$

*for $1 \le m < H_D$.*

(iii) *Each function $f_k$ is continuous, and satisfies a Lipschitz condition, on $D$.*

*Then the following are true.*

(a) *For $(\hat{z}_1, \dots, \hat{z}_q) \in D$, the system of differential equations*

$$\frac{dz_k}{dt} = f_k(z_1, \dots, z_q), \qquad k = 1, \dots, q, \tag{1}$$

*has a unique solution in $D$ for $z_k : \mathbb{R} \to \mathbb{R}$ passing through*

$$z_k(0) = \hat{z}_k, \qquad k = 1, \dots, q,$$

*and which extends to points arbitrarily close to the boundary of $D$.*

(b) *Let* $\lambda > \lambda_1$ *with* $\lambda = o(1)$ *and let* $\eta(\beta,\lambda) = \frac{\beta}{\lambda}\exp(-\frac{n\lambda^3}{\beta^3})$. *For a sufficiently large constant* $C$, *with probability* $1 - O(\eta(\beta,\lambda))$,

$$X_{n,m,k} = nz_k(m/n) + O(\lambda n)$$

*uniformly for* $0 \le m \le \sigma n \le m_n$ *and for each* $k$, *where* $z_k(t)$ *is the solution in* (a) *with* $\hat{z}_k = n^{-1}X_{n,0,k}$, *and* $\sigma = \sigma(n)$ *is the supremum of those* $m$ *to which the solution can be extended before reaching within* $L_\infty$*-distance* $C\lambda$ *of the boundary of* $D$.

We can now state our main theorem, which is based on Theorem 1. The multivariate normal distribution with mean 0 and covariance matrix $\Sigma$ is denoted by $\mathcal{N}(0,\Sigma)$.

THEOREM 2. *Assume that the conditions of Theorem 1 are satisfied, with* $\beta = o(n^{1/12-\varepsilon})$ *for some* $\varepsilon > 0$ *and* $\lambda_1 = o(n^{-1/2})$. *Furthermore, assume that the functions* $f_k$ *are differentiable, and that each partial derivative of* $f_k$ *is continuous, on* $D$. *Let* $z_1(t),\dots,z_q(t)$ *be the functions obtained in* (b) *of Theorem 1. Let* $g_{ij}:\mathbb{R}^q \to \mathbb{R}$, $1 \le i,j \le q$, *be functions, and assume that the following conditions hold.*

(ii′) *For some function* $\lambda_2 = \lambda_2(n) = o(1)$ *and all* $i,j$ *with* $1 \le i,j \le q$,

$$|\mathbb{E}[\Delta X_{n,m,i}\Delta X_{n,m,j} \mid \mathcal{F}_{m-1}] - g_{ij}(n^{-1}X_{n,m,1},\dots,n^{-1}X_{n,m,q})| \le \lambda_2$$

*for* $1 \le m < H_D$.

(iii′) *Each function* $g_{ij}$ *is continuous and satisfies a Lipschitz condition on* $D$.

*Then there is a continuous matrix-valued function* $\Sigma:\mathbb{R} \to \mathbb{R}^{q\times q}$ *such that*

$$\frac{\mathbf{X}_{n,m} - n\mathbf{z}(m/n)}{\sqrt{n}} \xrightarrow{d} \mathcal{N}(0,\Sigma(m/n)),$$

*where* $\mathbf{z}(t) = [z_1(t),\dots,z_q(t)]'$, *for* $0 \le m \le \sigma n$.

The proof of Theorem 2 in Section 3 also describes the procedure for calculating the matrix $\Sigma(t)$.

**2. A central limit theorem for near-martingales.** Our proof of Theorem 2 will be based on a central limit theorem for multidimensional martingales. Let $\{\mathbf{S}_{n,m}; m = 0,1,\dots,m_n\}$ be an array of random $q$-dimensional vectors with $\mathbf{S}_{n,0} = 0$. We denote the $k$th entry in $\mathbf{S}_{n,m}$ by $S_{n,m,k}$ and let as before $\Delta\mathbf{S}_{n,m} = \mathbf{S}_{n,m} - \mathbf{S}_{n,m-1}$. This theorem is the multidimensional version of Corollary 2.6 in [3].

THEOREM 3. *Let $\mathbf{S}_{n,m}$ be an array as above, and let $\Sigma = \{\sigma_{ij}\}_{i,j}$ be a $q \times q$-matrix. Assume that the following conditions are satisfied.*

(i) $\max_m \|\Delta \mathbf{S}_{n,m}\|$ *has uniformly bounded second moment.*
(ii) $\max_m \|\Delta \mathbf{S}_{n,m}\| \xrightarrow{p} 0$.
(iii) *For all* $1 \le i, j \le q$, $\sum_{m=1}^{m_n} \Delta S_{n,m,i} \Delta S_{n,m,j} \xrightarrow{p} \sigma_{ij}$.
(iv) $\sum_{m=1}^{m_n} \mathbb{E}[\Delta \mathbf{S}_{n,m} \mid \mathcal{F}_{m-1}] \xrightarrow{p} 0$.
(v) $\sum_{m=1}^{m_n} \mathbb{E}[\Delta \mathbf{S}_{n,m} \mid \mathcal{F}_{m-1}]^2 \xrightarrow{p} 0$.

*Then* $\mathbf{S}_{n,m_n} \xrightarrow{d} \mathcal{N}(0, \Sigma)$.

PROOF. Corollary 2.6 of McLeish [3] asserts that the theorem is true when $q = 1$ and $\sigma_{11} = 1$. It follows easily that the theorem also holds for arbitrary $\sigma_{11}$ in the univariate case.

Assume that $q > 1$, and let $\mathbf{a} = [a_1, \dots, a_q]' \in \mathbb{R}^q$ be an arbitrary $q$-dimensional vector. Let $R_{n,m} = \sum_{k=1}^q a_k S_{n,m,k}$. Since

$$\Delta R_{n,m} = \sum_{k=1}^q a_k \Delta S_{n,m,k},$$

it is easy to see that (i), (ii) and (iv) are satisfied for $R_{n,m}$. Assumption (v) means that $\sum_m \mathbb{E}[\Delta S_{n,m,i} \mid \mathcal{F}_{m-1}] \mathbb{E}[\Delta S_{n,m,j} \mid \mathcal{F}_{m-1}] \xrightarrow{p} 0$ for all $1 \le i, j \le q$. Hence,

$$\sum_m \mathbb{E}[\Delta R_{n,m} \mid \mathcal{F}_{m-1}]^2 = \sum_{i,j} a_i a_j \sum_m \mathbb{E}[\Delta S_{n,m,i} \mid \mathcal{F}_{m-1}] \mathbb{E}[\Delta S_{n,m,j} \mid \mathcal{F}_{m-1}]$$

tends to 0 in probability, so (v) holds also for $R_{n,m}$. Finally, we have

$$\begin{aligned} \sum_m (\Delta R_{n,m})^2 &= \sum_m \left( \sum_{k=1}^q a_k \Delta S_{n,m,k} \right)^2 = \sum_m \sum_{1 \le i,j \le q} a_i a_j \Delta S_{n,m,i} \Delta S_{n,m,j} \\ &= \sum_{1 \le i,j \le q} a_i a_j \sum_m \Delta S_{n,m,i} \Delta S_{n,m,j} \xrightarrow{p} \sum_{1 \le i,j \le q} a_i a_j \sigma_{ij}, \end{aligned}$$

so (iii) is satisfied for $R_{n,m_n}$ with parameter $\mathbf{a}' \Sigma \mathbf{a}$. Hence, by the univariate version of the theorem, $R_{n,m_n} \xrightarrow{d} \mathcal{N}(0, \mathbf{a}' \Sigma \mathbf{a})$. Since this holds for all vectors $\mathbf{a} \in \mathbb{R}^q$, it follows that $\mathbf{S}_{n,m_n} \xrightarrow{d} \mathcal{N}(0, \Sigma)$ (see, *e.g.*, Definition 3.2.5 in [5]). $\square$

**3. Proof of main theorem.** This section is devoted to the proof of Theorem 2. We are given a sequence of random $q$-dimensional vectors $\mathbf{X}_{n,m}$ and functions $f_k$ with $1 \le k \le q$ and $g_{ij}$ with $1 \le i, j \le q$ such that the conditions of Theorem 2 are satisfied. We will generally suppress $n$ in the subscript, so we write $\mathbf{X}_m$ for $\mathbf{X}_{n,m}$ and so on. It follows from the assumptions of Theorem 2 that we can choose a function $\lambda = o(n^{-1/4})$ such that $\lambda > \beta n^{-1/3+\varepsilon}$ for some $\varepsilon > 0$. Thus, according to Theorem 1, there are functions $\alpha_1(t), \ldots, \alpha_q(t)$ such that

$$X_{m,k} = n\alpha_k(m/n) + o(n^{3/4}) \tag{2}$$

with probability $1 - O(e^{-n^{\varepsilon}})$. Let $\mathcal{E}$ be the event that (2) holds for $1 \leq k \leq q$ and $0 \leq m < m_n$. Then $\mathbb{P}[\overline{\mathcal{E}}] = O(ne^{-n^{\varepsilon}})$. It is sufficient to prove that the conclusion of the theorem holds conditioned on $\mathcal{E}$. Indeed, let $\boldsymbol{\alpha}(t) = [\alpha_1(t), \ldots, \alpha_q(t)]'$ and $\mathbf{W}_m = n^{-1/2}(\mathbf{X}_m - n\boldsymbol{\alpha}(m/n))$. We have for an arbitrary bounded continuous function $\gamma$,

$$\begin{aligned}\mathbb{E}[\gamma(\mathbf{W}_m)] &= \mathbb{P}[\mathcal{E}]\mathbb{E}[\gamma(\mathbf{W}_m) \mid \mathcal{E}] + \mathbb{P}[\overline{\mathcal{E}}]\mathbb{E}[\gamma(\mathbf{W}_m) \mid \overline{\mathcal{E}}] \\ &= \mathbb{E}[\gamma(\mathbf{W}_m) \mid \mathcal{E}] + O(ne^{-n^{\varepsilon}}).\end{aligned}$$

Thus, if $\mathbf{W}_m$ tends to a normal distribution conditioned on $\mathcal{E}$, it also tends to a normal distribution when not conditioned on anything. In the following, we therefore assume that $\mathcal{E}$ holds.

Let $\mathbf{F}: \mathbb{R}^q \to \mathbb{R}^q$ be the vector-valued function whose $k$th component is $f_k$; that is,

$$\mathbf{F}(z_1, \ldots, z_q) = \begin{bmatrix} f_1(z_1, \ldots, z_q) \\ \vdots \\ f_q(z_1, \ldots, z_q) \end{bmatrix}.$$

By the assumption of Theorem 2, $\lambda_1 = o(n^{-1/2})$, so condition (ii) of Theorem 1 implies that

$$\mathbb{E}[\Delta\mathbf{X}_m \mid \mathcal{F}_{m-1}] = \mathbf{F}(n^{-1}\mathbf{X}_{m-1}) + o(n^{-1/2}). \tag{3}$$

We write $\boldsymbol{\alpha}_m = \boldsymbol{\alpha}(m/n)$. If we let $t = m/n$, then Taylor's theorem implies that

$$\begin{aligned}\boldsymbol{\alpha}_{m+1} &= \boldsymbol{\alpha}(t + n^{-1}) = \boldsymbol{\alpha}(t) + n^{-1}\frac{d\boldsymbol{\alpha}(t)}{dt} + O(n^{-2}) \\ &\overset{(1)}{=} \boldsymbol{\alpha}_m + n^{-1}\mathbf{F}(\boldsymbol{\alpha}_m) + O(n^{-2}),\end{aligned}$$

so

$$n\Delta\boldsymbol{\alpha}_m = \mathbf{F}(\boldsymbol{\alpha}_{m-1}) + O(n^{-1}), \tag{4}$$

analogous to (3). The Jacobian matrix of $\mathbf{F}$ is

$$J(\mathbf{z}) = \begin{bmatrix} \dfrac{\partial f_1}{\partial z_1} & \cdots & \dfrac{\partial f_1}{\partial z_q} \\ \vdots & \ddots & \vdots \\ \dfrac{\partial f_q}{\partial z_1} & \cdots & \dfrac{\partial f_q}{\partial z_q} \end{bmatrix}.$$

From calculus, we know that if $\mathbf{a}, \mathbf{y} \in \mathbb{R}^q$, then

$$\mathbf{F}(\mathbf{a}+\mathbf{y}) - \mathbf{F}(\mathbf{a}) = J(\mathbf{a})\mathbf{y} + O(\|\mathbf{y}\|^2) \tag{5}$$

as $\mathbf{y} \to 0$. We now let $\mathbf{Y}_m = \mathbf{X}_m - n\boldsymbol{\alpha}_m$ be the centered version of $\mathbf{X}_m$. By (2),

$$\mathbf{Y}_m = o(n^{3/4}). \tag{6}$$

Thus,

$$\begin{aligned} \mathbf{F}(n^{-1}\mathbf{X}_m) - \mathbf{F}(\boldsymbol{\alpha}_m) &= \mathbf{F}(\boldsymbol{\alpha}_m + n^{-1}\mathbf{Y}_m) - \mathbf{F}(\boldsymbol{\alpha}_m) \\ &\overset{(5)}{=} J(\boldsymbol{\alpha}_m)n^{-1}\mathbf{Y}_m + o(n^{-1/2}), \end{aligned} \tag{7}$$

so

$$\begin{aligned} \mathbb{E}[\Delta\mathbf{Y}_m] &= \mathbb{E}[\Delta\mathbf{X}_m] - n\mathbb{E}[\Delta\boldsymbol{\alpha}_m] \\ &\overset{(3,4)}{=} \mathbf{F}(n^{-1}\mathbf{X}_{m-1}) - \mathbf{F}(\boldsymbol{\alpha}_{m-1}) + o(n^{-1/2}) \\ &\overset{(7)}{=} J(\boldsymbol{\alpha}_{m-1})n^{-1}\mathbf{Y}_{m-1} + o(n^{-1/2}) \overset{(6)}{=} o(n^{-1/4}). \end{aligned}$$

Thus, $\mathbb{E}[\Delta\mathbf{Y}_m]$ tends to 0; however, the bound we have obtained is not strong enough to apply Theorem 3 directly to $\mathbf{Y}_m$. We will instead consider a transformation $\mathbf{Z}_m = T_m\mathbf{Y}_m$, where $T_m$ is a $q \times q$-matrix chosen so that $\mathbb{E}[\Delta\mathbf{Z}_m] = o(n^{-1/2})$ and $n^{-1}\sum_m Z_{m,i}Z_{m,j} \xrightarrow{p} \xi_{ij}(t)$ for some functions $\xi_{ij}(t)$. Then we will apply Theorem 3 to the array $n^{-1/2}\mathbf{Z}_m$, showing that it converges to a multivariate normal distribution. The normality of $\mathbf{X}_m$ will then be inferred from the normality of $\mathbf{Z}_m$.

For ease of notation, we write $A(t) = J(\boldsymbol{\alpha}(t))$. Note that $A(t)$ is a continuous matrix-valued function. Next, we define $T(t)$ to be the $q \times q$-matrix satisfying the differential equation

$$\frac{d}{dt}T(t) = -T(t)A(t), \qquad T(0) = I. \tag{8}$$

Lemma 1. *There is an open interval $(t_1, t_2)$ containing $[0, \sigma]$ such that there is a unique solution to the differential equation (8) on $(t_1, t_2)$, which furthermore satisfies a Lipschitz condition on $(t_1, t_2)$. If $T(t)$ satisfies (8), then $T(t)$ is invertible for all $t \in (t_1, t_2)$. Furthermore, let $T_m = T(m/n)$, $A_m = A(m/n)$ and*

$$U_m = I - n^{-1}A_m. \tag{9}$$

*Then*

$$T_{m+1} = T_mU_m + O(n^{-2}) \tag{10}$$

*for $0 \le m \le \sigma n - 1$.*

PROOF. Let $A(t) = \{a_{ij}(t)\}_{ij}$. If $\boldsymbol{\tau}_i(t) = [\tau_{i1}(t), \ldots, \tau_{iq}(t)]$ is the $i$th row of $T(t)$, then it is a solution of the system of linear homogenous differential equations

$$\frac{d}{dt}\tau_{ij}(t) = -\sum_{k=1}^{q} a_{kj}(t)\tau_{ij}(t), \qquad \tau_{ij}(0) = \delta_{ij}, \tag{11}$$

which can also be written $\boldsymbol{\tau}_i(t) = -\boldsymbol{\tau}_i(t)A(t)$, $\boldsymbol{\tau}_i(0) = \mathbf{e}_i$. Thus, every $\boldsymbol{\tau}_i(t)$ is actually a solution to the same system of linear differential equations; only the boundary condition is different.

Let $t_1' = \inf\{t : \boldsymbol{\alpha}(t') \in D \text{ for } t < t' < 0\}$ and $t_2' = \sup\{t : \boldsymbol{\alpha}(t') \in D \text{ for } 0 < t' < t\}$, and choose $t_1, t_2$ such that $t_1' < t_1 < 0 < \sigma < t_2 < t_2'$. By assumption, $A(t)$ is continuous on $(t_1', t_2')$. Hence, according to Theorem 12, Chapter 2 of Hurewicz [1], there is a unique solution to (11) on $(t_1', t_2')$. Moreover, by Theorem 2, Chapter 3 of [1], the solutions $\boldsymbol{\tau}_1(t), \ldots, \boldsymbol{\tau}_q(t)$ are linearly independent for all $t \in (t_1', t_2')$ if they are linearly independent for some $t \in (t_1', t_2')$. Thus, since $T(0) = I$ is invertible, $T(t)$ is invertible for all $t \in (t_1', t_2')$.

Since $A(t)$ and $T(t)$ are continuous on $(t_1', t_2')$, they are bounded on $(t_1, t_2)$. Thus, by (8), $\frac{d}{dt}T(t)$ is bounded on $(t_1, t_2)$, and so $T(t)$ satisfies a Lipschitz condition on $(t_1, t_2)$.

Finally we obtain by Taylor's theorem that

$$\begin{aligned} T_{m+1} &= T(t + n^{-1}) = T(t) + n^{-1}\frac{d}{dt}T(t) + O(n^{-2}) \\ &\overset{(8)}{=} T_m - n^{-1}T_m A_m + O(n^{-2}) \overset{(9)}{=} T_m U_m + O(n^{-2}). \end{aligned}$$

□

The matrices $A(t)$ and $T(t)$ do not depend on $n$, so we have $A(t), T(t) = O(1)$. As indicated, we now define $\mathbf{Z}_m = T_m \mathbf{Y}_m$. The next two lemmas show that $\mathbf{Z}_m$ has the properties required in order to apply Theorem 3 to the array $n^{-1/2}\mathbf{Z}_m$.

LEMMA 2. *For all* $m$,

$$\Delta \mathbf{Z}_m = O(\beta) \qquad \textit{a.s.}, \tag{12}$$

*and*

$$\mathbb{E}[\Delta \mathbf{Z}_m \mid \mathcal{F}_{m-1}] = o(n^{-1/2}). \tag{13}$$

PROOF. We have

$$\begin{aligned} \Delta \mathbf{Z}_m &= T_m \mathbf{Y}_m - T_{m-1}\mathbf{Y}_{m-1} \\ &\overset{(10)}{=} (T_{m-1}U_{m-1} + O(n^{-2}))\mathbf{Y}_m - T_{m-1}\mathbf{Y}_{m-1} \\ &\overset{(6)}{=} T_{m-1}(U_{m-1}\mathbf{Y}_m - \mathbf{Y}_{m-1}) + o(n^{-1}). \end{aligned} \tag{14}$$

By (4), $n\Delta\boldsymbol{\alpha}_m = O(1)$, so

$$\begin{aligned}\|\Delta\mathbf{Z}_m\| &= O(\|U_{m-1}\mathbf{Y}_m - \mathbf{Y}_{m-1}\|) \overset{(9)}{=} O(\|\Delta\mathbf{Y}_m\|) + O(n^{-1}\|A_{m-1}\mathbf{Y}_m\|) \\ &\le \|\Delta\mathbf{X}_m\| + \|n\Delta\boldsymbol{\alpha}_m\| + o(n^{-1/4}) = O(\beta) + O(1),\end{aligned}$$

implying (12). Then we consider (13), and first show that the conditional expectation of the term inside the parentheses in (14) is small. We have

$$\begin{aligned}&\mathbb{E}[U_{m-1}\mathbf{Y}_m - \mathbf{Y}_{m-1} \mid \mathcal{F}_{m-1}] \\ &\quad = U_{m-1}(\mathbb{E}[\mathbf{X}_m \mid \mathcal{F}_{m-1}] - n\boldsymbol{\alpha}_m) - \mathbf{X}_{m-1} + n\boldsymbol{\alpha}_{m-1} \\ &\quad = (I - n^{-1}A_{m-1})(\mathbf{X}_{m-1} + \mathbf{F}(n^{-1}\mathbf{X}_{m-1}) - n\boldsymbol{\alpha}_{m-1} - \mathbf{F}(\boldsymbol{\alpha}_{m-1})) \\ &\qquad - \mathbf{X}_{m-1} + n\boldsymbol{\alpha}_{m-1} + o(n^{-1/2}) \\ &\quad = -n^{-1}A_{m-1}(\mathbf{X}_{m-1} - n\boldsymbol{\alpha}_{m-1}) + \mathbf{F}(n^{-1}\mathbf{X}_{m-1}) - \mathbf{F}(\boldsymbol{\alpha}_{m-1}) + o(n^{-1/2}) \\ &\quad \overset{(7)}{=} -n^{-1}A_{m-1}\mathbf{Y}_{m-1} + n^{-1}A_{m-1}\mathbf{Y}_{m-1} + o(n^{-1/2}) \\ &\quad = o(n^{-1/2}),\end{aligned}$$

where we for the second equality have used (3), (4) and (9). Thus,

$$\begin{aligned}\mathbb{E}[\Delta\mathbf{Z}_m \mid \mathcal{F}_{m-1}] &\overset{(14)}{=} T_{m-1}\mathbb{E}[U_{m-1}\mathbf{Y}_m - \mathbf{Y}_{m-1} \mid \mathcal{F}_{m-1}] + o(n^{-1}) \\ &= T_{m-1}o(n^{-1/2}) + o(n^{-1}) = o(n^{-1/2}).\end{aligned}$$ □

We now turn to the quadratic variation.

LEMMA 3. *For all $m$,*

$$(\Delta\mathbf{Z}_m)^2 = O(\beta^2) \qquad a.s. \tag{15}$$

*Moreover, for $1 \le i, j \le q$, there is a function $\xi_{ij} : \mathbb{R} \to \mathbb{R}$, such that*

$$n^{-1}\sum_{k=1}^{m} \Delta Z_{k,i}\Delta Z_{k,j} \xrightarrow{p} \xi_{ij}(m/n). \tag{16}$$

PROOF. We have by (14) that

$$(\Delta\mathbf{Z}_m)^2 = T_{m-1}(U_{m-1}\mathbf{Y}_m - \mathbf{Y}_{m-1})^2 T'_{m-1} + o(1), \tag{17}$$

and by (9) that

$$(U_{m-1}\mathbf{Y}_m - \mathbf{Y}_{m-1})^2 = (\Delta\mathbf{Y}_m)^2 + o(1).$$

Since $n\Delta\boldsymbol{\alpha}_m = O(1)$, it follows from condition (i) of Theorem 1 that

$$\begin{aligned}(\Delta\mathbf{Y}_m)^2 &= (\Delta\mathbf{X}_m - n\Delta\boldsymbol{\alpha}_m)(\Delta\mathbf{X}_m' - n\Delta\boldsymbol{\alpha}_m') \\ &= (\Delta\mathbf{X}_m)^2 - \Delta\mathbf{X}_m n\Delta\boldsymbol{\alpha}_m' - n\Delta\boldsymbol{\alpha}_m\Delta\mathbf{X}_m' + n^2(\Delta\boldsymbol{\alpha}_m)^2 \\ &= O(\beta^2).\end{aligned} \tag{18}$$

This implies (15). To show (16), we take the conditional expectation and get

$$\begin{aligned}\mathbb{E}[(\Delta\mathbf{Y}_m)^2 \mid \mathcal{F}_{m-1}] &\overset{(18)}{=} \mathbb{E}[(\Delta\mathbf{X}_m)^2 \mid \mathcal{F}_{m-1}] - \mathbb{E}[\Delta\mathbf{X}_m \mid \mathcal{F}_{m-1}]n\Delta\boldsymbol{\alpha}_m' \\ &\quad - n\Delta\boldsymbol{\alpha}_m\mathbb{E}[\Delta\mathbf{X}_m' \mid \mathcal{F}_{m-1}] + n^2\Delta\boldsymbol{\alpha}_m^2 \\ &\overset{(3,4)}{=} \mathbb{E}[(\Delta\mathbf{X}_m)^2 \mid \mathcal{F}_{m-1}] - \mathbf{F}(n^{-1}\mathbf{X}_{m-1})\mathbf{F}(\boldsymbol{\alpha}_{m-1})' \\ &\quad - \mathbf{F}(\boldsymbol{\alpha}_{m-1})\mathbf{F}(n^{-1}\mathbf{X}_{m-1})' + \mathbf{F}(\boldsymbol{\alpha}_{m-1})^2 + o(1) \\ &\overset{(5)}{=} \mathbb{E}[(\Delta\mathbf{X}_m)^2 \mid \mathcal{F}_{m-1}] - \mathbf{F}(\boldsymbol{\alpha}_{m-1})^2 + o(1).\end{aligned}$$

Thus, by (17),

$$\begin{aligned}\mathbb{E}[(\Delta\mathbf{Z}_m)^2 \mid \mathcal{F}_{m-1}] &= T_{m-1}\mathbb{E}[(U_m\mathbf{Y}_m - \mathbf{Y}_{m-1})^2 \mid \mathcal{F}_{m-1}]T_{m-1}' + o(1) \\ &= T_{m-1}(\mathbb{E}[(\Delta\mathbf{X}_m)^2 \mid \mathcal{F}_{m-1}] - \mathbf{F}(\boldsymbol{\alpha}(t))^2)T_{m-1}' + o(1).\end{aligned} \tag{19}$$

Let $\mathbf{G}:\mathbb{R}^q \to \mathbb{R}^{q\times q}$ be the matrix-valued function such that

$$\mathbf{G}(z_1,\dots,z_q) = \{g_{ij}(z_1,\dots,z_q)\}_{i,j}.$$

Condition (ii$'$) of Theorem 2 can then be expressed as

$$\mathbb{E}[(\Delta\mathbf{X}_m)^2 \mid \mathcal{F}_{m-1}] = \mathbf{G}(n^{-1}X_{m,1},\dots,n^{-1}X_{m,q}) + o(1). \tag{20}$$

For $1 \le i,j \le q$, let $\zeta_{m,i,j} = \sum_{k=1}^m \Delta Z_{k,i}\Delta Z_{k,j}$, and let $Q_m = \{\zeta_{m,i,j}\}_{i,j} = \sum_{k=1}^m (\Delta\mathbf{Z}_k)^2$. Using (19) and (20), we find that if $t = m/n$, then

$$\begin{aligned}\mathbb{E}[Q_m - Q_{m-1} \mid \mathcal{F}_{m-1}] &= \mathbb{E}[(\Delta\mathbf{Z}_m)^2 \mid \mathcal{F}_{m-1}] \\ &= T(t)(\mathbf{G}(n^{-1}\mathbf{X}_{m-1}) - \mathbf{F}(\boldsymbol{\alpha}(t))^2)T(t)' + o(1).\end{aligned} \tag{21}$$

For $1 \le i,j \le q$, let $h_{ij}:\mathbb{R}^{q+1} \to \mathbb{R}$ be the functions such that

$$T(t)(\mathbf{G}(z_1,\dots,z_q) - \mathbf{F}(\boldsymbol{\alpha}(t))^2)T(t)' = \{h_{ij}(t,z_1,\dots,z_q)\}_{i,j}. \tag{22}$$

Then it follows from (21) that for $1 \le i,j \le q$,

$$|\mathbb{E}[\Delta\zeta_{m,i,j} \mid \mathcal{F}_{m-1}] - h_{ij}(m/n, X_{m,1}/n,\dots,X_{m,q}/n)| \le \lambda_3, \tag{23}$$

for some function $\lambda_3 = \lambda_3(n) = o(1)$.

Let $V_m$ be a random variable such that $V_m = m$ a.s. We will now apply Theorem 1 to the random variables in the set $\{V_m\} \cup \{X_{m,k}\}_k \cup \{\zeta_{m,i,j}\}_{i,j}$.

Since $\Delta V_m = 1$, the conditions of Theorem 1 are clearly satisfied by $V_m$. Moreover, we already know by assumption that they are satisfied by $X_{m,k}$ and $f_k$. Thus, we only have to check that they are also satisfied by $\zeta_{m,i,j}$ and $h_{ij}$.

By (15), $|\Delta\zeta_{m,i,j}| \le \|(\Delta\mathbf{Z}_m)^2\| = O(\beta^2)$, so condition (i) is satisfied. Condition (ii) is satisfied because of (23). To see that condition (iii) is satisfied, we have to show that the functions $h_{ij}$ are continuous and satisfy a Lipschitz condition on some area in $\mathbb{R}^{q+1}$.

Let $t_1$ and $t_2$ be as in Lemma 1. Let

$$D' = \{(t, z_1, \ldots, z_q) : t_1 < t < t_2, (z_1, \ldots, z_q) \in D\}.$$

Let us consider $\mathbf{F}$, $\mathbf{G}$ and $T$ as functions from $\mathbb{R}^{q+1}$ to $\mathbb{R}$, such that if $t \in \mathbb{R}$ and $\mathbf{z} \in \mathbb{R}^q$, then $\mathbf{F}(t,\mathbf{z}) = \mathbf{F}(\mathbf{z})$, $\mathbf{G}(t,\mathbf{z}) = \mathbf{G}(\mathbf{z})$ and $T(t,\mathbf{z}) = T(t)$. Since $D'$ is bounded, the product of two Lipschitz continuous functions on $D'$ is itself Lipschitz continuous on $D'$. By Lemma 1, $T(t,\mathbf{z})$ satisfies a Lipschitz condition, and by the assumptions, $\mathbf{F}(t,\mathbf{z})$ and $\mathbf{G}(t,\mathbf{z})$ do so as well. It then follows from the definition of $h_{ij}$ in (22) that $h_{ij}$ satisfies a Lipschitz condition on $D'$.

Let

$$\xi_{ij}(t) = \int h_{ij}(t)\,dt, \qquad \xi_{ij}(0) = 0.$$

Since $\beta^2 = o(n^{1/6})$, we can choose a function $\lambda' = o(1)$ such that $\lambda' > \lambda_3$ and $\eta(\beta^2, \lambda') = o(1)$. Then Theorem 1 implies that

$$\zeta_{m,i,j} = n\xi_{ij}(m/n) + O(\lambda' n),$$

for $0 \le m \le \sigma n$ and $1 \le i, j \le q$, with probability $1 - o(1)$. Hence, (16) is proved. $\square$

Lemma 4. *Let $\mathbf{M}_m = n^{-1/2}\mathbf{Z}_m$ and let $\Xi(t) = \{\xi_{ij}(t)\}_{i,j}$. Then*

$$\mathbf{M}_m \xrightarrow{d} \mathcal{N}(0, \Xi(m/n)).$$

Proof. We will show that $\mathbf{M}_m$ satisfies the conditions of Theorem 3.

(i) By (15) in Lemma 3,

$$\|(\Delta\mathbf{M}_k)^2\| = n^{-1}\|(\Delta\mathbf{Z}_k)^2\| = O(\beta^2/n) = o(1),$$

so $\max_k \|\Delta\mathbf{M}_k\|$ has uniformly bounded second moment.

(ii) By (12) in Lemma 2,

$$\max_k \|\Delta\mathbf{M}_k\| = n^{-1/2} \max_k \|\Delta\mathbf{Z}_k\| = O(\beta/\sqrt{n}) = o(1).$$

(iii) By (16) in Lemma 3,

$$\sum_{k=1}^{m} \Delta M_{k,i} \Delta M_{k,j} = n^{-1} \sum_{k=1}^{m} \Delta Z_{k,i} \Delta Z_{k,j} \xrightarrow{p} \xi_{ij}(t).$$

(iv) By (13) in Lemma 2,

$$\begin{aligned}\sum_{k=1}^{m} \mathbb{E}[\Delta \mathbf{M}_k \mid \mathcal{F}_{k-1}] &= n^{-1/2} \sum_{k=1}^{m} \mathbb{E}[\Delta \mathbf{Z}_k \mid \mathcal{F}_{k-1}] \\ &= n^{-1/2} m \cdot o(n^{-1/2}) = o(1).\end{aligned}$$

(v) Again by Lemma 2,

$$\begin{aligned}\sum_{k=1}^{m} \mathbb{E}[\Delta \mathbf{M}_k \mid \mathcal{F}_{k-1}]^2 &= n^{-1} \sum_{k=1}^{m} \mathbb{E}[\Delta \mathbf{Z}_k \mid \mathcal{F}_{k-1}]^2 \\ &= n^{-1} m \cdot o(n^{-1}) = o(1).\end{aligned}$$

The conclusion then follows from Theorem 3. □

Proof of Theorem 2. From Lemma 1, we know that $T(t)$ is invertible, so we can define $\Sigma(t) = T(t)^{-1} \Xi(t) (T(t)^{-1})'$. We then conclude from Lemma 4 that

$$\frac{\mathbf{X}_{n,m} - n\boldsymbol{\alpha}(m/n)}{\sqrt{n}} \xrightarrow{d} \mathcal{N}(0, \Sigma(m/n)).$$

□

**4. Random graph processes with restricted degrees.** For a positive integer $d$, the random $d$-process is a random graph process defined as follows. Begin with an empty graph on $n$ vertices. Every step in the process consists of choosing two distinct vertices in the graph uniformly at random, and adding an edge between them if and only if the vertices are not adjacent and both of them have degree at most $d-1$. The process ends when the graph no longer contains a pair of nonadjacent vertices, both of which have degree smaller than $d$. It was proved in [4] that the graph process asymptotically almost surely (i.e., with probability tending to 1 as $n \to \infty$, abbreviated a.a.s.) produces a graph where at most one vertex has degree $d-1$ while all other vertices have degree $d$. If $dn$ is even, the final graph is a.a.s. $d$-regular.

This process was used in [6] to illustrate the usage of the differential equation method. Here, we show that the present central limit theorem also can be applied to the process. Let $G_m$ be the graph after $m$ edges have been added, and let $V_{m,k}$ be the random variable denoting the number of vertices of degree $k$ in $G_m$. We follow the argument in [6] and note that

$$\mathbb{E}[\Delta V_{m,k} \mid \mathcal{F}_{m-1}] = \frac{2\delta_{k>0} V_{m-1,k-1} - 2\delta_{k<d} V_{m-1,k}}{n - V_{m-1,d}} + o(1).$$

Moreover, $|\Delta V_{m,k}| \le 2$ always, and the domain $D$ is chosen as $-\varepsilon < z_i < 1+\varepsilon$ for $0 \le k < d$ and $\varepsilon < z_d < 1-\varepsilon$ for some $\varepsilon > 0$. All the conditions of Theorem 1 are therefore satisfied, and it follows that there are functions $\gamma_0(t), \dots, \gamma_d(t)$ such that a.a.s.

$$V_{m,k} = \gamma_k(m/n)n + o(n)$$

for $0 \le k \le d$. In order to apply Theorem 2, we note that

$$\begin{aligned}&\mathbb{E}[\Delta V_{m,i}\Delta V_{m,j} \mid \mathcal{F}_{m-1}]\\ &\quad = \sum_{k=0}^{d-1}\sum_{l=0}^{d-1} \frac{V_{m-1,k}}{n - V_{m-1,d}} \frac{V_{m-1,l}}{n - V_{m-1,d}} (\delta_{k,i-1} - \delta_{k,i})(\delta_{l,j-1} - \delta_{l,j}) + o(1),\end{aligned}$$

so condition (ii′) of Theorem 2 holds for the functions

$$g_{ij}(z_1, \dots, z_q) = \sum_{k=0}^{d-1}\sum_{l=0}^{d-1} \frac{z_k z_l}{(1-z_d)^2} (\delta_{k,i-1} - \delta_{k,i})(\delta_{l,j-1} - \delta_{l,j}).$$

We choose $D$ to be the same as earlier, and note that the functions $g_{ij}$ satisfy a Lipschitz condition on $D$. Theorem 2 then implies the following theorem.

THEOREM 4. *Let* $\mathbf{V}_m = [V_{m,0}, \dots, V_{m,d}]'$ *and* $\boldsymbol{\gamma}(t) = [\gamma_0(t), \dots, \gamma_d(t)]'$. *Let* $\delta > 0$ *and let* $m_\delta$ *be the smallest value for which* $\gamma_d(m_\delta/n) > 1 - \delta$. *There is a continuous matrix-valued function* $\Sigma(t)$ *such that*

$$\frac{\mathbf{V}_m - n\boldsymbol{\gamma}(m/n)}{\sqrt{n}} \xrightarrow{d} \mathcal{N}(0, \Sigma(m/n))$$

*for* $0 \le m \le m_\delta$.

**5. The minimum degree random graph process.** Our second application is the first phase of the minimum-degree graph process, first introduced in [7]. One complication in this case is that the graph process has a natural random stopping time, and we will show that the random variables under consideration also have a jointly normal distribution at the end of the process.

For a fixed $n$, the minimum degree graph process is a sequence of graphs $\{G_m^{\min}\}_{m\ge 0}$ which is constructed as follows. The initial graph $G_0^{\min}$ is an empty graph on $n$ vertices. For $m \ge 1$, let $v_m$ be a vertex chosen uniformly at random from the vertices of minimum degree in $G_{m-1}^{\min}$, and let $w_m$ be chosen uniformly at random from the vertices distinct from $v_m$. The graph $G_m^{\min}$ is obtained from $G_{m-1}^{\min}$ by adding to it the edge $(v_m, w_m)$. For simplicity, we will allow multi-edges; however, in the stages of the process we consider, there will a.a.s. be so few multi-edges that they make no significant difference to the calculations.

Let $H$ be such that $G_H^{\min}$ does not contain isolated vertices, while $G_{H-1}^{\min}$ contains at least one isolated vertex. In [7], it was proved that a.a.s. $H = hn + o(n)$, where $h = \ln 2$. In this paper, we will consider the graph only up to the point $H$; that is, we add edges at random until there are no isolated vertices left, and then we stop. Thus, $H$ is a stopping time of the process, and we consider the process $G_{H\wedge m}^{\min}$. It is easy to see that no cycle can be formed before time $H$, so $G_H^{\min}$ is a forest. In [7], it was furthermore proved that the number of vertices of any degree is sharply concentrated around the expectation. Instead of the vertex degrees, we will consider the order of the components in $G_m^{\min}$, and in particular in $G_H^{\min}$. For $k \geq 1$ and $m \geq 0$, let $C_{m,k}$ be the random variable denoting the number of components in $G_m^{\min}$ of order $k$, and let $C_k = C_{H,k}$. Let $\mathbf{C}_m = [C_{m,1}, \ldots, C_{m,q}]'$ and $\mathbf{C} = [C_1, \ldots, C_q]'$, where $q \geq 1$ is some fixed natural number. In [2], it was shown that $C_{m,k} = \beta_k(m/n)n + o(n)$ a.a.s., where

$$\beta_k(t) = \frac{1}{k}(1 - e^{-t})^{k-1}((k+1)e^{-t} - 1).$$

Let $\boldsymbol{\beta}(t) = [\beta_1(t), \ldots, \beta_q(t)]'$. We will prove the following theorem.

Theorem 5. *There is a continuous matrix-valued function $\Sigma(t)$ such that*

$$\frac{\mathbf{C}_m - n\boldsymbol{\beta}(m/n)}{\sqrt{n}} \xrightarrow{d} \mathcal{N}(0, \Sigma(m/n)) \tag{24}$$

*for $0 \leq \frac{m}{n} < h$.*

*Let $\boldsymbol{\mu} = \{\frac{k-1}{k2^k}\}_{k=1}^q$. Then there is matrix $\Sigma$ such that*

$$\frac{\mathbf{C} - n\boldsymbol{\mu}}{\sqrt{n}} \xrightarrow{d} \mathcal{N}(0, \Sigma). \tag{25}$$

Proof. Assume first that $\frac{m}{n} = t < h$, where $t$ is a constant. When a new edge $(v_m, w_m)$ is added, $v_m$ is by definition an isolated vertex, while $w_m$ can have any degree, and be in a component of any order. Let $V_m$ be the random variable denoting the order of the component containing $w_m$. Then

$$\Delta C_{m,k} = -\delta_{k1} - \delta_{k,V_m} + \delta_{k-1,V_m}. \tag{26}$$

The probability of choosing a vertex in a component of order $k$ is

$$\mathbb{P}[V_m = k] = \frac{kC_{m-1,k} - \delta_{k1}}{n-1},$$

so the expected change in the number of components of order $k$ is

$$\mathbb{E}[\Delta C_{m,k} \mid \mathcal{F}_{m-1}] = f_k(n^{-1}C_{m,1}, \ldots, n^{-1}C_{m,q}) + o(n^{-1/2}),$$

where

$$f_k(z_1, \dots, z_q) = -\delta_{k1} - kz_k + (k-1)z_{k-1}. \tag{27}$$

Furthermore, it is clear that $C_{m,k} \le n$ and (26) implies that

$$\Delta C_{m,k} \le 2 \tag{28}$$

for $m \ge 1$ and $1 \le k \le q$. The set $D$ can be chosen as $\varepsilon < z_1 < 1 + \varepsilon$ and $-\varepsilon < z_k < 1$ for $2 \le k \le q$ for any $\varepsilon > 0$. Then $f_k$ satisfy a Lipschitz condition on $D$. We obtain a system of differential equations of the form (1), and it can be shown that it has the solution $z_k = \beta_k(t)$, satisfying the boundary conditions $\beta_k(0) = \delta_{k1}$. Let $t_0 = h - \delta$. For every $\delta > 0$, we can choose $\varepsilon$ so small that the solution does not leave $D$ until $t > t_0$. It follows that a.a.s.

$$C_{m,k} = \beta_k(t)n + o(n)$$

for $1 \le k \le q$ and $0 \le t < h$, with $t$ fixed. This was already shown in [2].

In order to apply Theorem 2, we need an expression for the conditional expectation of $\Delta C_{m,i}\Delta C_{m,j}$. This is

$$\mathbb{E}[\Delta C_{m,i}\Delta C_{m,j} \mid \mathcal{F}_{m-1}] = g_{ij}(n^{-1}C_{m,1}, \dots, n^{-1}C_{m,q}) + o(1),$$

where

$$g_{ij}(z_1, \dots, z_q) = \sum_{k \ge 1} kz_k(-\delta_{i1} - \delta_{ki} + \delta_{k,i-1})(-\delta_{j1} - \delta_{jk} + \delta_{k,j-1}). \tag{29}$$

We let $D$ be as earlier. Then the conditions of Theorem 2 are satisfied, and we conclude that there is a matrix $\Sigma(t)$ such that (24) holds.

We then turn to (25). Let $m_0 = \lfloor t_0 n \rfloor$. By (28), $|C_{H,k} - C_{m_0,k}| \le 2\delta n$. By letting $\delta$ go to 0, we can conclude that

$$n^{-1}C_{H,k} \xrightarrow{p} \beta_k(h) = \frac{k-1}{k2^k}.$$

Unfortunately, we cannot obtain the matrix $\Sigma$ simply by setting $t = h$ in $\Sigma(t)$, since the stopping time $H$ is a random variable. Instead, we define new random variables $C^\circ_{m,1}, \dots, C^\circ_{m,q}$ which behave like $C_{m,k}$ up to $m = H$, but which we can analyze even after $H$. Let $V^\circ_m$ for $m \ge 1$ be defined as follows. If $C_{m-1,1} > 0$, let $V^\circ_m = V_m$. Otherwise, let $V^\circ_m = k$ with probability $\frac{kC^\circ_{m,k}}{n}$ for $2 \le k \le q$ and $q+1$ with probability $1 - \sum_{k=2}^q \frac{kC^\circ_{m,k}}{n}$. Let $C^\circ_{0,k} = n\delta_{k1}$ and define

$$C^\circ_{m,k} = C^\circ_{m-1,k} - \delta_{k1} - \delta_{k,V^\circ_m} + \delta_{k-1,V^\circ_m}.$$

Clearly, $C^\circ_{m,k} = C_{m,k}$ for $0 \le m \le H$. We observe that

$$\mathbb{E}[\Delta C^\circ_{m,k} \mid \mathcal{F}_{m-1}] = f_k(n^{-1}C^\circ_{m,1}, \dots, n^{-1}C^\circ_{m,q}) + o(n^{-1/2})$$

and

$$\mathbb{E}[\Delta C^\circ_{m,i}\Delta C^\circ_{m,j} \mid \mathcal{F}_{m-1}] = g_{ij}(n^{-1}C^\circ_{m,1},\dots,n^{-1}C^\circ_{m,q}) + o(1)$$

when $m \le hn + o(n)$, where $f_k$ and $g_{ij}$ are defined by (27) and (29), respectively. Let $D^\circ \subset \mathbb{R}^q$ be the defined by $-\varepsilon < z_k < 1+\varepsilon$ for $1 \le k \le q$. Then the requirements of Theorem 2 are satisfied, and we can conclude that

$$\frac{\mathbf{C}^\circ_{hn} - n\boldsymbol{\mu}}{\sqrt{n}} \xrightarrow{d} \mathcal{N}(0,\Sigma(h)).$$

However, we are interested in the distribution of $\mathbf{C}_H = \mathbf{C}^\circ_H$, so we have to find the difference between $\mathbf{C}^\circ_H$ and $\mathbf{C}^\circ_{hn}$. For $1 \le k \le q$, let $W^\circ_k$ be random variables such that $[W^\circ_1,\dots,W^\circ_q]' \sim \mathcal{N}(0,\Sigma(h))$. Thus, for example, $n^{-1/2}C^\circ_{hn,1} \xrightarrow{d} W^\circ_1$. Let $\eta_n = \frac{H-hn}{\sqrt{n}}$. Then writing $m' = m - hn$,

$$\begin{aligned}\eta_n > \frac{m-hn}{\sqrt{n}} &\Leftrightarrow H > m \Leftrightarrow C'_{m,1} > 0\\ &\Leftrightarrow W^\circ_1\left(\frac{m}{n}\right) > -\frac{2e^{-m/n}-1}{\sqrt{n}} = -\frac{e^{-m'/n}-1}{\sqrt{n}}\\ &\Leftrightarrow W^\circ_1\left(\frac{m}{n}\right) > (1+o(1))\frac{m-hn}{\sqrt{n}},\end{aligned}$$

so $\eta_n \xrightarrow{d} W^\circ_1$.

When $m = hn + o(n)$, $\mathbb{P}[V^\circ_m = k] = \frac{k-1}{2^k} + o(1)$ for $1 \le k \le q$. Let

$$B_k = \operatorname{sgn}(H - hn)\sum_{H\wedge hn}^{H\vee hn} \delta_{V_m,k}.$$

Thus, $|B_k|$ is the number of times a vertex in a component of order $k$ is chosen between $H$ and $hn$. Then $\mathbb{E}[B_k] = (H-hn)(\frac{k-1}{2^{k-1}} + o(1))$ and one easily sees that $\mathbb{E}[|B_k|(|B_k|-1)] = (1+o(1))\mathbb{E}[|B_k|]^2$, so by Chebyshev's inequality, $B_k \sim (H-hn)\frac{k-1}{2^{k-1}}$. Hence, $n^{-1/2}B_k \xrightarrow{d} \frac{k-1}{2^{k-1}}\eta_n$, and we can conclude that

$$\frac{C_k - n(k-1)/(k2^k)}{\sqrt{n}} \xrightarrow{d} W_k,$$

where $W_k = W^\circ_k + W^\circ_1\frac{k-1}{2^{k-1}}$ for $1 \le k \le q$. Since $\{W_1,\dots,W_q\}$ are linear combinations of $\{W^\circ_1,\dots,W^\circ_q\}$, they are jointly normal random variables. □

When the functions $f_k$ are linear, as in this section, it becomes easier to calculate $\Sigma(t)$ explicitly, than in the nonlinear case. The matrix $A$ is then a constant matrix and the solution of (8) is $T(t) = e^{-tA}$, where the matrix exponential is defined as $e^{tA} = \sum_{i\ge 0}\frac{(tA)^i}{i!}$.

In our example, the Jacobian matrix is $J = \{j\delta_{i,j+1} - j\delta_{ij}\}_{ij}$, and $A = J$, so

$$T(t) = \{\delta_{j\leq i}(-1)^{i+j}e^{jt}(e^t - 1)^{i-j}\}_{ij}.$$

The covariance matrix is then given by

$$\Sigma(t) = e^{tA}\int e^{-tA}(\mathbf{G}(\boldsymbol{\beta}(t)) - \mathbf{F}(\boldsymbol{\beta}(t))^2)e^{-tA'}\,dt e^{tA'},$$

which can be used as a general formula for the linear case.

## REFERENCES

[1] Hurewicz, W. (1958). *Lectures on Ordinary Differential Equations*. The Technology Press of the Massachusetts Institute of Technology, Cambridge, MA. MR0090703

[2] Kang, M. and Seierstad, T. G. (2007). Phase transition of the minimum degree random multigraph process. *Random Structures Algorithms* **31** 330–353. MR2352179

[3] McLeish, D. L. (1974). Dependent central limit theorems and invariance principles. *Ann. Probab.* **2** 620–628. MR0358933

[4] Ruciński, A. and Wormald, N. C. (1992). Random graph processes with degree restrictions. *Combin. Probab. Comput.* **1** 169–180. MR1179247

[5] Tong, Y. L. (1990). *The Multivariate Normal Distribution*. Springer, New York. MR1029032

[6] Wormald, N. C. (1995). Differential equations for random processes and random graphs. *Ann. Appl. Probab.* **5** 1217–1235. MR1384372

[7] Wormald, N. C. (1999). The differential equation method for random graph processes and greedy algorithms. In *Lectures on Approximation and Randomized Algorithms*. (M. Karoński and H. J. Prömel, eds.) 75–152. PWN, Warsaw.

Department of Biostatistics
Institute of Basic Medical Sciences
University of Oslo
Postboks 1122 Blindern
N-0317 Oslo
Norway
E-mail: taralgs@medisin.uio.no
URL: http://folk.uio.no/taralgs